\documentclass[a4paper,DIV=12,10pt]{scrartcl}

\usepackage[utf8]{inputenc}
\usepackage{amsfonts,amssymb,amsmath}
\usepackage{mathtools}

\usepackage{epsfig}
\usepackage{MnSymbol}
\usepackage{booktabs}
\usepackage{cite}
\usepackage{float}
\usepackage{hyperref}
\usepackage{siunitx}

\usepackage{wrapfig}

\usepackage{graphics}
\usepackage{color}

\numberwithin{equation}{section}
\numberwithin{table}{section}
\numberwithin{figure}{section}

\newcommand{\ud}{\,\mathrm{d}}
\newcommand{\R}{\mathbb{R}}
\newcommand{\N}{\mathbb{N}}

\newcommand*{\ldblbrace}{\{\mskip-5mu\{}
\newcommand*{\rdblbrace}{\}\mskip-5mu\}}

\newcommand{\disc}{{\operatorname{disc}}}
\renewcommand{\vec}{\boldsymbol}

\numberwithin{equation}{section}
\numberwithin{figure}{section}
\numberwithin{table}{section}

\newtheorem{defi}{Definition}[section]
\newtheorem{prob}[defi]{Problem}
\newtheorem{rem}[defi]{Remark}

\usepackage{cleveref}

\usepackage{caption}
\usepackage{subcaption}

\begin{document}
	
	\title{Benchmark computations of dynamic poroelasticity}
	
	\author{
		Mathias Anselmann$^\dag$, Markus Bause$^\dag$\thanks{bause@hsu-hh.de (corresponding author)}\;, Nils Margenberg$^\dag$, Pavel Shamko$^\dag$\\
		{\small $^\dag$ Helmut Schmidt University, Faculty of
			Mechanical and Civil Engineering, Holstenhofweg 85,}\\ 
		{\small 22043 Hamburg, Germany}
	}
	\date{}
	
	\maketitle

	\begin{abstract}
We present benchmark computations of dynamic poroelasticity modeling fluid flow in deformable porous media by a coupled hyperbolic-parabolic system of partial differential equations. A challenging benchmark setting and goal quantities of physical interest for this problem are proposed. Computations performed by space-time finite element approximations with continuous and discontinuous discretizations of the time variable are summarized. By this work we intend to stimulate comparative studies by other research groups for the evaluation of dynamic poroelasticity solver regarding the accuracy of discretization techniques, the efficiency and robustness of iterative methods for the linear systems and the arrangement of the model equations in terms of their variables (two-field or multi-field formulations). 
	\end{abstract}

\section{Introduction and mathematical model}

In this work we present benchmark computations for families of space-time finite element approximations (cf., e.g., \cite{ABMS23,BAKR22,HST13,HST11,KB14,KM04}) to the coupled hyperbolic-parabolic problem
\begin{subequations}
	\label{Eq:HPS}
	\begin{alignat}{7}
		\label{Eq:HPS_1}
		\rho \partial_t^2 \vec u - \nabla \cdot (\vec C \vec \varepsilon (\vec u)) + \alpha \vec \nabla p & = \rho \vec f\,, \qquad 
		&
		c_0\partial_t p + \alpha \nabla \cdot \partial_t \vec u  -  \nabla \cdot (\vec K \nabla p)  & = g\,, &&\quad  \text{in } \; \Omega \times (0,T]\,,\\[1ex]
		\label{Eq:HPS_3}
	\vec u (0) = \vec u_0\,, \quad \partial_t \vec u (0) & = \vec u_1\,, & p(0) & = p_0\,, && \quad \text{in } \; \Omega\times \{0\} \,,\\[1ex]
		\label{Eq:HPS_4}
		 \vec u = \vec u_D\,, \quad \text{on } \; \Gamma_{\vec u}^{D}  \times  (&0,T]\,, & \quad 
-(\vec C\vec \varepsilon(\vec u) - \alpha p\vec E) \vec n  & = \vec t_N\,,  && \quad \text{on } \; \Gamma_{\vec u}^{{N}} \times (0,T]\,,\\[1ex]
\label{Eq:HPS_6}
 p = p_D\,, \quad \text{on } \; \Gamma_p^{{D}} \times  (&0,T]\,, & 
- \vec K \nabla p \cdot \vec n  & = p_N\,, && \quad \text{on } \; \Gamma_p^{{N}} \times (0,T]\,.
	\end{alignat}
\end{subequations}
Componentwise or directional boundary conditions for $\vec u$, given by 
\begin{equation}
	\label{Eq:HPS_DBC}
		\vec u \cdot \vec n = 0 \qquad \text{and} \qquad  (\vec \sigma(\vec u)\vec n) \cdot \vec t_i  = 0\,, \quad \text{for}\;\; i=1,\ldots, d-1\,,  \quad \text{on}\;\; \Gamma^d_{\vec u} \times (0,T]\,,
\end{equation}
are applied further for the sake of physical realism; cf.~Fig.~\ref{fig:L_shaped_domain}. In \eqref{Eq:HPS}, $\Omega \subset \R^d$, with $d\in \{2,3\}$, is an open bounded Lipschitz domain with outer unit normal vector $\vec n$ to the boundary $\partial \Omega$ and $T>0$ is the final time point. For \eqref{Eq:HPS_4} and  \eqref{Eq:HPS_6}, we let $\partial \Omega = {\Gamma_{\vec u}^{D}} \cup {\Gamma_{\vec u}^{N}}$ and $\partial \Omega = {\Gamma_{p}^{D}}\cup {\Gamma_{p}^{N}}$ with closed portions $\Gamma_{\vec u}^{D}$ and $\Gamma_{p}^{D}$ of non-zero measure. In \eqref{Eq:HPS_DBC}, we denote by $\vec t_i$, for $i=1,\ldots,d-1$, the unit basis vectors of the tangent space at $\vec x \in \Gamma^d_{\vec u}$. For \eqref{Eq:HPS_DBC}, the decomposition $\partial \Omega = {\Gamma_{\vec u}^{N}} \cup \Gamma^d_{\vec u}$ of the boundary of a L-shaped domain $\Omega$ that is used for our computations is illustrated in Fig.~\ref{fig:L_shaped_domain}. The unknowns in \eqref{Eq:HPS} are the vector-valued variable $\vec u$ and the scalar function $p$. The quantity $\vec \varepsilon (\vec u):= (\nabla \vec u + (\nabla \vec u)^\top)/2$ denotes the symmetrized gradient and $\vec E\in \R^{d,d}$ is the identity matrix. For brevity, the positive quantities $\rho>0$, $\alpha>0$ and $c_0 >0$ as well as the tensors $\vec C$ and $\vec K$ are assumed to be constant in space and time. The tensors $\vec C$ and $\vec K$ are supposed to be symmetric and positive definite,
\begin{equation*}
	\label{Eq:PosDef}
		\exists k_0>0 \; \forall  \vec \xi = \vec \xi^\top \in \R^{d,d}:  \; \sum_{i,j,k,l=1}^d \xi_{ij} C_{ijkl} \xi_{kl} \geq k_0 \sum_{j,k=1}^d |\xi_{jk}|^2\,,\quad 
		\exists k_1>0 \; \forall \vec \xi \in \R^d:  \; \sum_{i,j,=1}^d \xi_{i} K_{ij} \xi_{j} \geq k_1 \sum_{i=1}^d |\xi_{i}|^2\,.
\end{equation*}
Under these assumptions, well-posedness of \eqref{Eq:HPS} is ensured. This has been shown by different mathematical techniques and for several combinations of boundary conditions in, e.g., \cite{JR18,S89,STW22}. In \cite{ABMS23}, well-posedness of a fully discrete space-time finite element approximation is proved for the boundary conditions in \eqref{Eq:HPS_DBC}, applied to the L-shaped domain $\Omega$ of Fig.~\ref{fig:L_shaped_domain}.
 
Important and classical applications of the system \eqref{Eq:HPS} arise in poroelasticity and thermoelasticity; cf., e.g., \cite{B41,B55,B72} and \cite{C72,JR18,L86}. Recently, generalizations of the system \eqref{Eq:HPS} to soft materials have strongly attracted researchers' interest in biomedicine; cf., e.g., \cite{BKNR22,CADQ22} and the references therein. In neurophysiology, such generalizations are used to model, simulate and elucidate circulatory diseases, such as ischaemic stroke, or also Alzheimer's disease. In poroelasticity, Eqs.~\eqref{Eq:HPS} are referred to as the dynamic Biot model. The system \eqref{Eq:HPS} is used to describe flow of a slightly compressible viscous fluid through a deformable porous matrix. The small deformations of the matrix are described by the Navier equations of linear elasticity, and the diffusive fluid flow is described by Duhamel’s equation. The unknowns are the effective solid phase displacement $\vec u$ and the effective fluid pressure $p$. The quantity $\vec \varepsilon (\vec u)$ is the strain tensor. Further, $\rho$ is the effective mass density, $\vec C$ is Gassmann’s fourth order effective elasticity tensor, $\alpha$ is Biot’s pressure-storage coupling tensor, $c_0$ is the specific storage coefficient and $\vec K$ is the permeability. In thermoelasticity, $p$ denotes the temperature, $c_0$ is the specific heat of the medium, and $\vec K$ is the conductivity. Then, the quantity $\alpha \nabla p$ arises from the thermal stress in the structure, and $\alpha \nabla \cdot \partial_t \vec u$ models the internal heating due to the dilation rate. 

\section{Space-time finite element approximation}
\label{Sec:Disc}

We rewrite \eqref{Eq:HPS} as a first-order in time system by introducing the new variable $\vec v:= \partial_t \vec u$. Then, we recover \eqref{Eq:HPS_1} as
\begin{equation}
	\label{Eq:HPS_11}
		\partial_t \vec u - \vec v = \vec 0\,, \quad 
		\rho \partial_t \vec v - \nabla \cdot (\vec C \vec \varepsilon (\vec u)) + \alpha \vec \nabla p = \rho \vec f\,, 
		\qquad 		c_0\partial_t p + \alpha \nabla \cdot \vec v  - \nabla \cdot (\vec K \vec \nabla p)  = g\,,
\end{equation}
along with the initial and boundary conditions \eqref{Eq:HPS_3} to \eqref{Eq:HPS_6}. We benchmark the application of space-time finite element methods to \eqref{Eq:HPS_11}, with continuous and discontinuous Galerkin methods for the discretization of the time variable and inf-sup stable pairs of finite elements for the approximation of the space variables. To introduce the scheme we need notation.

For the time discretization, we decompose $I:=(0,T]$ into $N$ subintervals $I_n=(t_{n-1},t_n]$, $n=1,\ldots,N$, where $0=t_0<t_1< \cdots < t_{N-1} < t_N = T$ such that $I=\bigcup_{n=1}^N I_n$. We put $\tau := \max_{n=1,\ldots, N} \tau_n$ with $\tau_n = t_n-t_{n-1}$. The set $\mathcal{M}_\tau := \{I_1,\ldots,
I_N\}$ of time intervals is called the time mesh. For a Banach space $B$, any $k\in
\N_0$ and $\mathbb P_k(I_n;B) := \big\{w_\tau \,: \,  I_n \to B \,, \; w_\tau(t) = \sum_{j=0}^k 
	W^j t^j \; \forall t\in I_n\,, \; W^j \in B\; \forall j \big\}$ we define the space of piecewise polynomial functions in time with values in $B$ by 
\begin{equation}
	\label{Eq:DefYk}
	Y_\tau^{k} (B) := \left\{w_\tau: \overline I \rightarrow B  \mid w_\tau{}_{|I_n} \in
	\mathbb P_{k}(I_n;B)\; \forall I_n\in \mathcal{M}_\tau,\, w_\tau(0)\in B \right\}\subset L^2(I;B)\,.
\end{equation}
For any function $w: \overline I\to B$ that is piecewise sufficiently smooth with respect to the time mesh $\mathcal{M}_{\tau}$, for instance for $w\in Y^k_\tau (B)$, we define the right-hand sided and left-hand sided limit at a mesh point $t_n$ by $w^+(t_n) := \lim_{t\to t_n+0} w(t)$ for $n<N$ and $w^-(t_n) := \lim_{t\to t_n-0} w(t)$ for $n>0$.  Further, for a Banach space $B$ and any $k\in \N$ we define 
\begin{equation}
	\label{Eq:DefXk}
	X_\tau^{k} (B) := \left\{w_\tau \in C(\overline I;B) \mid w_\tau{}_{|I_n} \in
	\mathbb P_{k}(I_n;B)\; \forall I_n\in \mathcal{M}_\tau \right\}\,.
\end{equation}
For time integration, it is natural to use the right-sided $(k+1)$-point Gau{ss}--Radau (GR) quadrature formula  in a discontinuous Galerkin approach and  the $(k+1)$-point Gau{ss}--Lobatto (GL) quadrature formula  in a continuous one. On $I_n$, they read as 
\begin{equation}
	\label{Eq:GF}
	Q_n^{\text{q}}(w) := \frac{\tau_n}{2}\sum_{\mu=1}^{k+1} \hat 
	\omega_\mu^{\operatorname{q}} w(t_{n,\mu}^{\operatorname{q}} ) \approx \int_{I_n} w(t) \ud t \,, \quad \text{for}\; \text{q}\in \{\operatorname{GR},\operatorname{GL}\}\,,
\end{equation}
where $t_{n,\mu}^{\operatorname{q}}=T_n(\hat t_{\mu}^{\operatorname{q}})$, for $\mu = 1,\ldots,k+1$, are the quadrature  points on $I_n$ and $\hat \omega_\mu^{\operatorname{q}}$ the corresponding weights of the respective quadrature formula. Here, $T_n(\hat t):=(t_{n-1}+t_n)/2 + (\tau_n/2)\hat t$ is the affine transformation from $\hat I = [-1,1]$ to $I_n$ and $\hat t_{\mu}^{\operatorname{q}}$ are the quadrature points on $\hat I$. Formula \eqref{Eq:GF} is exact for all $w\in \mathbb P_{2k} (I_n;\R)$ if $\text{q}=\operatorname{GR}$, and for all $w\in \mathbb P_{2k-1} (I_n;\R)$ if $\text{q}=\operatorname{GL}$ . 

For the space discretization, let $\mathcal{T}_h=\{K\}$  be the quasi-uniform decomposition of $\Omega$ into (open) quadrilaterals or hexahedrals, with mesh size $h>0$. These element types are chosen for our implementation (cf.\ Sec.~\ref{Sec:NumExp1} and \ref{Sec:NumExp2}) that uses the deal.II library \cite{Aetal21}.  The finite element spaces used for approximating the unknowns $\vec u$, $\vec v$ and $p$ of \eqref{Eq:HPS_11} in space are of the form 
\begin{subequations}
	\label{Def:VhQh}
	\begin{alignat}{2}
	\label{Def:Vh}
		\vec V_{h}  & := \{\vec v_h \in C(\overline \Omega )^d \mid \vec v_{h}{}_{|K}\in 
		{\vec V(K)} \;\; \text{for all}\; K \in \mathcal{T}_h\}\,, \\[1ex]
			\label{Def:Qh}
		Q_{h}  & := \{\vec q_h \in L^2(\Omega)\mid  \vec q_{h}{}_{|K}\in {Q(K)} \;\; \text{for all}\; K \in \mathcal{T}_h\}\,.
		\end{alignat}
\end{subequations}
For the local spaces $\vec V(K)$ and $Q(K)$ we employ mapped versions of the inf-sup stable pair $\mathbb Q_r^d/\mathbb P_{r-1}^\disc$, for $r\geq 2$, of finite element spaces; cf.\ \cite{J16}. The pair $\mathbb Q_r^d/\mathbb P_{r-1}^\disc$ with a discontinuous approximation of $p$ in the broken polynomial space $Q_h$ has proved excellent accuracy and stability properties for higher-order approximations of mixed (or saddle point) systems like the Navier--Stokes equations and the applicability of geometric multigrid preconditioner  for the algebraic systems; cf., e.g., \cite{AB23,AB22}.  

For $\vec w_h, \vec \chi_h\in \vec V_h$ and $q_h, \psi_h \in Q_h$ we define the bilinear forms 
\begin{subequations}
	\label{Def:ACB}
	\begin{alignat}{2}
		\hspace*{-3ex} A_\gamma (\vec w_{h},\vec \chi_h)	& := \langle \vec C \vec \varepsilon (\vec w_h),\vec \varepsilon (\vec \chi_h)\rangle - \langle \vec C\vec \varepsilon(\vec w_h) \vec n, \vec \chi_h\rangle_{\Gamma^D_{\vec u}}+ a_\gamma (\vec w_{h},\vec \chi_h)\,,
		\label{Def:A}
		\\[2ex]
		\hspace*{-3ex} C (\vec \chi_{h},q_h)   & := -\alpha \langle \nabla \cdot \vec \chi_h, q_h\rangle +  \alpha \langle \vec \chi_h \cdot \vec n , q_h \rangle_{\Gamma^D_{\vec u}} \,, 
		\label{Def:C}
		\\[2ex]
		\hspace*{-3ex} B_\gamma (q_{h},\psi_h) & :=\hspace*{-1ex}
			\begin{array}[t]{@{}l} \displaystyle \sum_{K\in \mathcal T_l} \hspace*{-0.5ex} \langle \vec K \nabla q_h, \nabla \psi_h \rangle_{K} - \hspace*{-1ex} \sum_{F\in \mathcal F_h} \hspace*{-1ex} \big (\langle \ldblbrace \vec K \nabla q_h \rdblbrace \cdot \vec n, \lsem \psi_h \rsem \rangle_{F}  + \langle \lsem q_h\rsem, \ldblbrace \vec K \nabla \psi_h \rdblbrace \cdot \vec n \rangle_{F}\big) + b_\gamma (q_{h},\psi_h)  \,,
			\end{array} 
        \label{Def:B}				
	\end{alignat}
\end{subequations}
where $a_\gamma (\cdot,\cdot)$ in \eqref{Def:A} is given by $
		a_\gamma (\vec w ,\vec \chi_h) := - \langle \vec w, \vec C\vec \varepsilon(\vec \chi_h) \vec n \rangle_{\Gamma^D_{\vec u}}  + {\gamma_a}{h_F^{-1}}  \langle \vec w, \vec \chi_h \rangle_{\Gamma^D_{\vec u}}$, 
for $\vec w \in \vec H^{1/2}(\Gamma^D_{\vec u})$, and $b_\gamma (\cdot,\cdot)$ in \eqref{Def:B} is defined by $b_\gamma (q_{h},\psi_h) = \sum_{F\in \mathcal F_h} {\gamma_b}{h_F^{-1}} \langle \lsem  q_h \rsem, \lsem \psi_h\rsem \rangle_F$. The form $B_\gamma$
yields a symmetric interior penalty discontinuous Galerkin discretization of the scalar variable $p$; cf., e.g.,  \cite[Sec.~4.2]{PE12}. As usual, the average $\ldblbrace \cdot \rdblbrace$  and jump $\lsem \cdot \rsem$ for a function $w$ of a broken space on an interior face $F$ between two elements $K^+$ and $K^-$, such that $F=\partial K^+ \cap \partial K^-$, are defined by $
	\ldblbrace w \rdblbrace :=  \frac{1}{2} (w^{+}+ w^{-})$ and $\lsem w \rsem : =  w^{+} -  w^{-} $. 
For boundary faces $F \subset \partial K \cap \partial \Omega$, we set $\ldblbrace w\rdblbrace:= w_{|K}$ and $\lsem w \rsem := w_{|K}$. The set of all faces (interior and boundary faces) on $\mathcal T_h$ is denoted by $\mathcal F_h$. In \eqref{Def:B}, the parameter $\gamma_b$ in $b_\gamma$ has to be chosen sufficiently large, such that the discrete coercivity of $B_\gamma$ on $Q_h$  is preserved. The local length $h_F$ is chosen {as $h_F = \ldblbrace h_F \rdblbrace := \frac{1}{2} (|K^+|_d + |K^-|_d)$} with Hausdorff measure $|\cdot |_d$; cf.~\cite[p.~125]{PE12}. For boundary faces we set $h_F := |K|_d$. In $a_\gamma (\cdot,\cdot)$ , the quantity $\gamma_a$ is the algorithmic parameter of the stabilization (or penalization) term in the Nitsche formulation \cite{B02,ABMS23} for incorporating Dirichlet boundary conditions in weak form which is applied here. To ensure well-posedness of the discrete systems, the parameter $\gamma_a$ has to be chosen  sufficiently large as well; cf.\ \cite[Appendix]{ABMS23}. {Based on our numerical experiments we choose the algorithmic parameter $\gamma_a$ and $\gamma_b$ as $\gamma_a = 5\cdot 10^4  \cdot r \cdot (r+1)$ and $\gamma_b =  \frac{1}{2} \cdot r \cdot (r-1)$, where $r$ is the polynomial degree of $\vec V_h$ in \eqref{Def:Vh}. Finally, for $\vec f \in \vec H^{-1}(\Omega)$, $\vec u_D \in \vec H^{1/2}(\Gamma^D_{\vec u})$, $\vec t_N \in \vec H^{-1/2}(\Gamma^N_{\vec u})$ and $g\in H^{-1}(\Omega)$, $p_D \in H^{1/2}(\Gamma^D_{p})$, $p_N\in H^{-1/2}(\Gamma^N_{p})$ we put 
\begin{subequations}
	\label{Def:LG}	
	\begin{alignat}{2}
		\nonumber 
		F_\gamma (\vec \chi_h)  := \langle \vec f , \vec \chi_h\rangle - \langle \vec t_N,\vec \chi_h \rangle_{\Gamma^N_{\vec u}} + a_\gamma (\vec u_D ,\vec \chi_h)\,,  \qquad 
		G_\gamma (\psi_h)   := \langle g,\psi_h \rangle -  \sum_{F\in \mathcal F_h^{D,\vec u}} \alpha \langle \vec v_D \cdot \vec n , \psi_h \rangle_{F} \\[1ex] 
		\nonumber 
		 -  \sum_{F\in \mathcal F_h^{D,p}}\langle p_D, \ldblbrace \vec K \nabla \psi_h \rdblbrace \cdot \vec n \rangle_{F} + \sum_{F\in \mathcal F_h^{D,p}} {\gamma_b}{h_F^{-1}} \langle p_D, \lsem \psi_h\rsem \rangle_F - \sum_{F\subset \mathcal F_h^{N,p}} \langle p_N, \lsem \psi_h \rsem \rangle_F\,.
	\end{alignat}
\end{subequations}
Here, we denote by $\mathcal F_h^{D,p}\subset \mathcal F_h$ and $\mathcal F_h^{N,p}\subset \mathcal F_h$ the set of all element faces on the boundary parts $\Gamma_p^D$ and  $\Gamma_p^N$, respectively; cf.\ \eqref{Eq:HPS}. The second of the terms on the right-hand side of $G_\gamma(\cdot,\cdot)$ with $\vec v_D = \partial_t \vec u_D$, is added to ensure consistency of the form $C(\cdot,\cdot)$ in \eqref{Def:C}; cf.\ \cite[p.\ 8]{ABMS23} for details.

We use a temporal test basis that is supported on the subintervals $I_n$. Then, a time marching process is obtained. In that, we assume that the trajectories $\vec u_{ \tau,h}$, $\vec v_{ \tau,h}$ and $p_{ \tau,h}$ have been computed before for all $t\in [0,t_{n-1}]$, starting with approximations  $\vec u_{\tau,h}(t_0) :=\vec u_{0,h}$, $\vec v_{\tau,h}(t_0) :=\vec u_{1,h}$ and $p_{\tau,h}(t_0) := p_{0,h}$ of the initial values $\vec u_0$, $\vec u_1$ and $p_{0}$. Then, we consider solving the following local problems on $I_n$ of the discontinuous (dG($k$)) and continuous (cG($k$)) Galerkin approximation in time; cf.\ \cite{ABMS23,BAKR22}. 

\begin{prob}[$I_n$-problem for dG($k$)]
	\label{Prob:DG}
	Let $k\in \N_0$. For given $\vec u_{h}^{n-1}:= \vec u_{\tau,h}(t_{n-1})\in \vec V_h$, $\vec v_{h}^{n-1}:=  \vec v_{\tau,h}(t_{n-1})\in \vec V_h$,  and $p_{h}^{n-1}:= p_{\tau,h}(t_{n-1}) \in Q_h$ with  $\vec u_{\tau,h}(t_0) :=\vec u_{0,h}$, $\vec v_{\tau,h}(t_0) :=\vec u_{1,h}$ and $p_{\tau,h}(t_0) := p_{0,h}$, find $(\vec u_{\tau,h},\vec v_{\tau,h},p_{\tau,h}) \in \mathbb P_k (I_n;\vec V_h) \times \mathbb P_k (I_n;\vec V_h) \times \mathbb P_k (I_n;Q_h)$ such that 	for all $(\vec \phi_{\tau,h},\vec \chi_{\tau,h},\psi_{\tau,h})\in  \mathbb P_k (I_n;\vec V_h) \times \mathbb P_k (I_n;\vec V_h) \times \mathbb P_k (I_n;Q_h)$, 
	\begin{subequations}
		\label{Eq:DG}
		\begin{alignat}{2}
			\label{Eq:DG1}
			&\begin{aligned}
				& Q_n^{\operatorname{GR}} \big(\langle \partial_t \vec u_{\tau,h} , \vec \phi_{\tau,h} \rangle  - \langle \vec v_{\tau,h} , \vec \phi_{\tau,h} \rangle \big) + \langle \vec u^+_{\tau,h}(t_{n-1}), \vec \phi_{\tau,h}^+(t_{n-1})\rangle =  \langle \vec u_{h}^{n-1}, \vec \phi_{\tau,h}^+(t_{n-1})\rangle\,, \\[0.5ex]
			\end{aligned}\\
			\label{Eq:DG2}
			&\begin{aligned}
				& Q_n^{\operatorname{GR}} \Big(\langle \rho \partial_t \vec v_{\tau,h} , \vec \chi_{\tau,h} \rangle + A_\gamma(\vec u_{\tau,h}, \vec \chi_{\tau,h} ) + C(\vec \chi_{\tau,h},p_{\tau,h})\Big) + \langle \rho \vec v^+_{\tau,h}(t_{n-1}),  \vec \chi_{\tau,h}^+(t_{n-1})\rangle\\[0.5ex]  
				& \quad = Q_n^{\operatorname{GR}} \Big(F_\gamma(\vec \chi_{\tau,h})\Big) + \langle \rho \vec v_{h}^{n-1},  \chi_{\tau,h}^+(t_{n-1})\rangle \,,\\[0.5ex]
			\end{aligned}\\
			\label{Eq:DG3}
			&\begin{aligned}
				& Q_n^{\operatorname{GR}} \Big(\langle c_0 \partial_t p_{\tau,h},\psi_{\tau,h} \rangle  - C(\vec v_{\tau,h},\psi_{\tau,h})+ B_\gamma (p_{\tau,h}, \psi_{\tau,h})\Big) + \langle c_0 p^+_{\tau,h}(t_{n-1}), \psi_{\tau,h}^+(t_{n-1})\rangle \\[0.5ex] 
				&  \quad = Q_n^{\operatorname{GR}} \Big( G_\gamma(\psi_{\tau,h})\Big) + \langle c_0 p_{h}^{n-1}, \psi_{\tau,h}^+(t_{n-1})\rangle\,.
			\end{aligned}
		\end{alignat}
	\end{subequations}
\end{prob}

The trajectories defined by Problem \ref{Prob:DG}, for $n = 1,\ldots,N$, satisfy that $\vec u_{\tau,h},\vec v_{\tau,h}\in Y_\tau^k(\vec V_h)$ and $p_{\tau,h}\in Y_\tau^k(Q_h)$; cf.\ \eqref{Eq:DefYk}. Well-posedness of Problem \ref{Prob:DG} is ensured; cf.~\cite[Lem.\ 3.2]{ABMS23}.

\begin{prob}[$I_n$-problem for cG($k$)]
	\label{Prob:CG}
	Let $k\in \N$. For given $\vec u_{h}^{n-1}:= \vec u_{\tau,h}(t_{n-1})\in \vec V_h$, $\vec v_{h}^{n-1}:=  \vec v_{\tau,h}(t_{n-1})\in \vec V_h$,  and $p_{h}^{n-1}:= p_{\tau,h}(t_{n-1}) \in Q_h$ with  $\vec u_{\tau,h}(t_0) :=\vec u_{0,h}$, $\vec v_{\tau,h}(t_0) :=\vec u_{1,h}$ and $p_{\tau,h}(t_0) := p_{0,h}$, find $(\vec u_{\tau,h},\vec v_{\tau,h},p_{\tau,h}) \in \mathbb P_k (I_n;\vec V_h) \times \mathbb P_k (I_n;\vec V_h) \times \mathbb P_k (I_n;Q_h)$ such that 
	\begin{equation}
	\label{Eq:CG0}
	\vec u_{\tau,h}^+(t_{n-1})=\vec u_{h}^{n-1}\,, \quad \vec v_{\tau,h}^+(t_{n-1})=\vec v_{h}^{n-1}\quad \text{and} \quad  p_{\tau,h}^+(t_{n-1})= p_{h}^{n-1}
	\end{equation}
	 and, for all $(\vec \phi_{\tau,h},\vec \chi_{\tau,h},\psi_{\tau,h})\in  \mathbb P_{k-1} (I_n;\vec V_h) \times \mathbb P_{k-1} (I_n;\vec V_h) \times \mathbb P_{k-1} (I_n;Q_h)$, 
	\begin{subequations}
		\label{Eq:CG}
		\begin{alignat}{2}
			\label{Eq:CG1}
			&\begin{aligned}
				& Q_n^{\operatorname{GL}} \big(\langle \partial_t \vec u_{\tau,h} , \vec \phi_{\tau,h} \rangle  - \langle \vec v_{\tau,h} , \vec \phi_{\tau,h} \rangle \big)  = 0 \,, \\[0.5ex]
			\end{aligned}\\
			\label{Eq:CG2}
			&\begin{aligned}
				& Q_n^{\operatorname{GL}} \Big(\langle \rho \partial_t \vec v_{\tau,h} , \vec \chi_{\tau,h} \rangle + A_\gamma(\vec u_{\tau,h}, \vec \chi_{\tau,h} ) + C(\vec \chi_{\tau,h},p_{\tau,h})\Big)  = Q_n^{\operatorname{GL}} \Big(F_\gamma(\vec \chi_{\tau,h})\Big) \,,\\[0.5ex]
			\end{aligned}\\
			\label{Eq:CG3}
			&\begin{aligned}
				& Q_n^{\operatorname{GL}} \Big(\langle c_0 \partial_t p_{\tau,h},\psi_{\tau,h} \rangle  - C(\vec v_{\tau,h},\psi_{\tau,h})+ B_\gamma (p_{\tau,h}, \psi_{\tau,h})\Big) = Q_n^{\operatorname{GL}} \Big( G_\gamma(\psi_{\tau,h})\Big)\,.
			\end{aligned}
		\end{alignat}
	\end{subequations}
	
\end{prob}

The trajectories defined by Problem \ref{Prob:CG}, for $n = 1,\ldots,N$, satisfy that $\vec u_{\tau,h},\vec v_{\tau,h}\in X_\tau^k(\vec V_h)$ and $p_{\tau,h}\in X_\tau^k(Q_h)$; cf.\ \eqref{Eq:DefXk}. Well-posedness of Problem \ref{Prob:CG} can be shown following the arguments \cite[Lem.\ 3.2]{ABMS23}. Precisely, the cG($k$) scheme of Problem~\ref{Prob:CG} represents a Galerkin--Petrov approach since trial and test spaces differ from each other.

\begin{rem}[Algebraic solver]
Problem~\ref{Prob:DG} and \ref{Prob:CG} lead to large linear algebraic systems with complex block structure, in particular for larger values of the piecewise polynomial order in time $k$. This puts a facet of complexity on their solution, in particular if three space dimensions are involved. To solve such type of block systems, we use GMRES iterations that are preconditioned by a V-cycle geometric multigrid method based on a local Vanka smoother; cf.\ \cite{AB23,ABMS23} for details.    
\end{rem}

\section{Numerical convergence test}
\label{Sec:NumExp1}
 
 \begin{table}[t!]
 	\centering
 	\begin{tabular}{l}
 		\begin{tabular}{cccccccc}
 			\toprule
 			{$\tau$} & {$h$} &
 			{ $\| \nabla (\vec u - \vec u_{\tau,h})  \|_{L^2(\vec L^2)} $ } & {EOC} &
 			{ $\| \vec v - \vec v_{\tau,h}  \|_{L^2(\vec L^2)} $ } & {EOC} &
 			{ $\| p-p_{\tau,h}  \|_{L^2(L^2)}  $ } & {EOC} 
 			\\
 			\cmidrule(r){1-2}
 			\cmidrule(lr){3-8}
 			$\tau_0/2^0$ & $h_0/2^0$ & 1.2138632264e-02 & {--} & 3.4963867086e-02  & {--} & 2.0325417612e-03 & {--} \\ 
 			$\tau_0/2^1$ & $h_0/2^1$ & 1.4699245816e-03 & 3.05 & 3.9349862997e-03  & 3.15 & 2.3314471675e-04 & 3.12 \\
 			$\tau_0/2^2$ & $h_0/2^2$ & 1.8238666739e-04 & 3.01 & 4.8313770355e-04  & 3.03 & 2.8891798035e-05 & 3.01 \\
 			$\tau_0/2^3$ & $h_0/2^3$ & 2.2707201873e-05 & 3.01 & 6.0087160884e-05  & 3.01 & 3.6067034147e-06 & 3.00 \\
 			$\tau_0/2^4$ & $h_0/2^4$ & 2.8305414233e-06 & 3.00 & 7.4900993341e-06  & 3.00 & 4.5021336614e-07 & 3.00 \\
 			\bottomrule
 		\end{tabular}
 	\end{tabular}
 	\caption{%
 		$L^2(L^2)$ errors and experimental orders of convergence (EOC) for \eqref{Eq:givensolution} with spatial degree $r=4$ for the spaces $\mathbb Q_r^2/\mathbb P_{r-1}^{\text{disc}}$  and temporal degree $k=2$ for the dG($k$) scheme of Problem \ref{Prob:DG}.
 	}
 	\label{Tab:1}
 \end{table}
 \nopagebreak[4]
 \begin{table}[h!]
 	\centering
 	\begin{tabular}{l}
 		\begin{tabular}{cccccccc}
 			\toprule
 			{$\tau$} & {$h$} &
 			{ $\| \nabla (\vec u - \vec u_{\tau,h})  \|_{L^2(\vec L^2)} $ } & {EOC} &
 			{ $\| \vec v - \vec v_{\tau,h}  \|_{L^2(\vec L^2)} $ } & {EOC} &
 			{ $\| p-p_{\tau,h}  \|_{L^2(L^2)}  $ } & {EOC} 
 			\\
 			\cmidrule(r){1-2}
 			\cmidrule(lr){3-8}
 			$\tau_0/2^0$ & $h_0/2^0$ & 9.8743046869e-04 & {--} & 3.5668679054e-03  & {--} & 4.3639594011e-04 & {--} \\ 
 			$\tau_0/2^1$ & $h_0/2^1$ & 5.9913786816e-05 & 4.04 & 1.5360551492e-04  & 4.54 & 2.5681365609e-05 & 4.09 \\
 			$\tau_0/2^2$ & $h_0/2^2$ & 3.7323826100e-06 & 4.00 & 9.0006987407e-06  & 4.09 & 1.5681328845e-06 & 4.03 \\
 			$\tau_0/2^3$ & $h_0/2^3$ & 2.3306835992e-07 & 4.00 & 5.5529338177e-07  & 4.02 & 9.7664800443e-08 & 4.01 \\
 			$\tau_0/2^4$ & $h_0/2^4$ & 1.4566772574e-08 & 4.00 & 3.4466105658e-08  & 4.01 & 6.1040500343e-09 & 4.00 \\
 			\bottomrule
 		\end{tabular}
 	\end{tabular}
 	\caption{%
 		$L^2(L^2)$ errors and experimental orders of convergence (EOC) for \eqref{Eq:givensolution} with spatial degree $r=4$ for the spaces $\mathbb Q_r^2/\mathbb P_{r-1}^{\text{disc}}$ and temporal degree $k=3$ for the cG($k$) scheme of Problem \ref{Prob:CG}.
 	}
 	\label{Tab:2}
 \end{table}

Firstly, we investigate the schemes proposed in Problem~\ref{Prob:DG} and \ref{Prob:DG} by a numerical convergence study. From the point of view of numerical costs for solving the algebraic counterparts of \eqref{Eq:DG} and \eqref{Eq:CG0}, \eqref{Eq:CG}, respectively, the dG($k$) member of the family of schemes in Problem~\ref{Prob:DG} can be compared with the cG($k$+1) scheme of Problem~\ref{Prob:CG}. If  in cG($k$) the Gauss--Lobatto quadrature points are used for building Lagrange interpolation in time on $\overline I_n$, the degrees of freedom at time $t_{n-1}$ are directly obtained  from the vector identities corresponding to  \eqref{Eq:CG0}. Using this, the algebraic system of cG($k$) can be condensed. Then, the dimension of the resulting algebraic system for cG($k$+1) coincides with the one obtained for dG($k$).

We study \eqref{Eq:HPS} for $\Omega=(0,1)^2$ and $I=(0,2]$ and the prescribed solution
\begin{equation}
	\label{Eq:givensolution}
	\boldsymbol u(\boldsymbol x, t) = \phi(\boldsymbol x, t) \boldsymbol E_2 \;\; \text{and}\;\;
	p(\boldsymbol x, t) = \phi(\boldsymbol x, t) \;\; \text{with}\;\; 
	\phi(\boldsymbol x, t) = \sin(\omega_1 t^2) \sin(\omega_2 x_1) \sin(\omega_2 x_2)
\end{equation}
and $\omega_1=\omega_2 = \pi$. We put $\rho=1.0$, $\alpha=0.9$, $c_0=0.01$ and $\boldsymbol K=\boldsymbol E_2$ with the identity $\vec E_2\in \R^{2,2}$. For the fourth order elasticity tensor $\boldsymbol C$, isotropic material properties with Young's modulus $E=100$ and Poisson's ratio $\nu=0.35$, corresponding to the Lam\'e parameters $\lambda = 86.4$ and $\mu = 37.0$, are chosen. For the space-time convergence test, the domain $\Omega$ is decomposed into a sequence of successively refined meshes of quadrilateral finite elements. The spatial and temporal mesh sizes are halved in each of the refinement steps. The step sizes of the coarsest space and time mesh are $h_0=1/(2\sqrt{2})$ and $\tau_0=0.1$. 

For the dG($k$) scheme of Problem~\ref{Prob:DG} we choose the polynomial degrees $k=2$ and $r=4$, such that discrete solutions $\vec u_{\tau,h}, \vec v_{\tau}\in Y_\tau^2(\vec V_h)$, $p_{\tau,h}\in Y_\tau^2(Q_h)$ with local spaces $\mathbb Q_4^2/ P_3^{\text{disc}}$ are obtained. For the cG($k$) scheme of Problem~\ref{Prob:CG} we choose the polynomial degrees $k=3$ and $r=4$, such that discrete solutions $\vec u_{\tau,h}, \vec v_{\tau}\in X_\tau^3(\vec V_h)$, $p_{\tau,h}\in X_\tau^3(Q_h)$ with local spaces $\mathbb Q_4^2/ P_3^{\text{disc}}$ are obtained. The calculated errors and corresponding experimental orders of convergence are summarized in Tab.~\ref{Tab:1} and \ref{Tab:2}, respectively. The error is measured in the quantities associated with the energy of the system \eqref{Eq:HPS}; cf.\ \cite[p.~15]{JR18} and \cite{BAKR22}. Table~\ref{Tab:2} nicely confirms the optimal rates of convergence with respect to the polynomial degrees in space and time for the cG($3$) scheme. The superiority of the cG($3$) scheme over the dG($2$) scheme is clearly observed.  We note that $r=3$ would have been sufficient for dG(2) to ensure third order convergence in space and time, $r=4$ was only chosen to equilibrate the costs for solving the algebraic systems and, thereby, compare either approaches fairly to each other.

\section{Benchmark computations }
\label{Sec:NumExp2}

Here we propose the two-dimensional case of our benchmark problem for dynamic poroelasticity. We intend to stimulate other research groups to use this benchmark for the evaluation of their schemes and implementations and contribute to its dissemination and, possibly, further improvement. The test problem is expected to enable comparative studies of different formulations of \eqref{Eq:HPS} in terms of unknowns (two-field versus multi-field arrangements) and to benchmark numerical approaches with respect to their accuracy and efficiency. 

We consider the L-shaped domain $\Omega\subset \R^3$ sketched in Fig.~\ref{fig:L_shaped_domain} along with the boundary conditions for $\vec u$ prescribed on the different parts of $\partial \Omega$. Beyond the boundary conditions \eqref{Eq:HPS_4}, we apply the (homogeneous) directional boundary conditions \eqref{Eq:HPS_DBC} on the portion $\Gamma_{\vec u}^d$ of $\partial \Omega$.  For their implementation in the forms of Subsec.~\ref{Sec:Disc} we refer to \cite[Subsec.~5.2]{ABMS23}. We aim to compute goal quantities of physical interest that are defined by 
\begin{equation}
	\label{Eq:DefGQ}
	G_{\vec{u}} = \int_{\Gamma_m} \vec{u} \cdot \vec{n} \, \ud o \qquad \text{and} \qquad 
	G_{p} = \int_{\Gamma_m} p \, \ud o \qquad \text{for}\;\; \Gamma_m:= \{0.75\}\times (0,0.5)\,.
\end{equation} 
On the left part of the upper boundary, i.e.\ for $\Gamma^N_{\vec u}:=(0,0.5)\times \{1\}$,  we impose the traction force 
\begin{equation*}
	\vec  t_N =  \begin{cases}
	-64x^2(16x-3)\sin(8\pi t)\,, & x \in [0,\tfrac{1}{8}]\,, \\[1ex]
	\frac{16}{27}(2x-1)^2(16x+1)\sin(8\pi t)\,, & x \in (\tfrac{1}{8},0.5]\,.
\end{cases}
\end{equation*}  
On the right boundary we put $\vec t_N=\vec 0$. For the variable $p$ we prescribe a homogeneous Dirichlet condition \eqref{Eq:HPS_6} on the left upper part of $\partial \Omega$, i.e., for $(x_1,x_2)\in \Gamma_p^D:=[0,0.5]\times \{1\}$. On $\Gamma_p^N:=\partial \Omega \backslash \Gamma_p^D$ we impose the  homogeneous Neumann condition $p_N=0$ in \eqref{Eq:HPS_6}. We put $\rho=1.0$, $\alpha=0.9$, $c_0=0.01$ and $\boldsymbol K=\boldsymbol E_2$ with the identity $\vec E_2\in \R^{2,2}$. For the elasticity tensor $\boldsymbol C$, isotropic material properties with Young's modulus $E=20000$ and Poisson's ratio $\nu=0.3$ are used. The final time is $T=8$.

\begin{figure}[t!]
	\centering
	\subcaptionbox{Test setting with boundary conditions for $\vec u$. \label{fig:L_setting}}
	[0.34\columnwidth]
	{\includegraphics[width=0.2\textwidth,keepaspectratio]{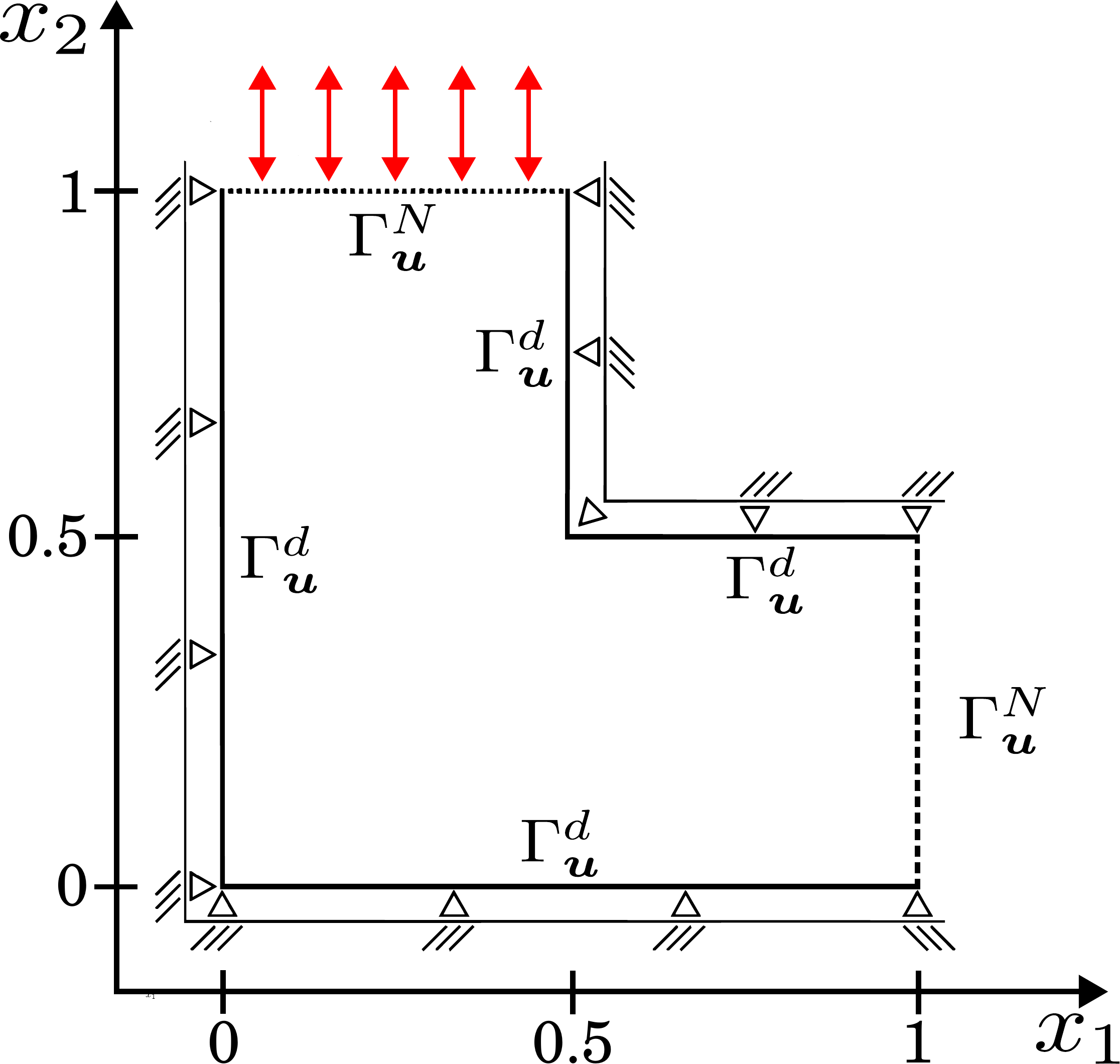}
	}
	\subcaptionbox{Modulus $|\vec u|$ of $\vec u$ at time t = 3.0. \label{fig:L_solution}}
	[0.25\columnwidth]
	{\includegraphics[width=0.17\textwidth,keepaspectratio]{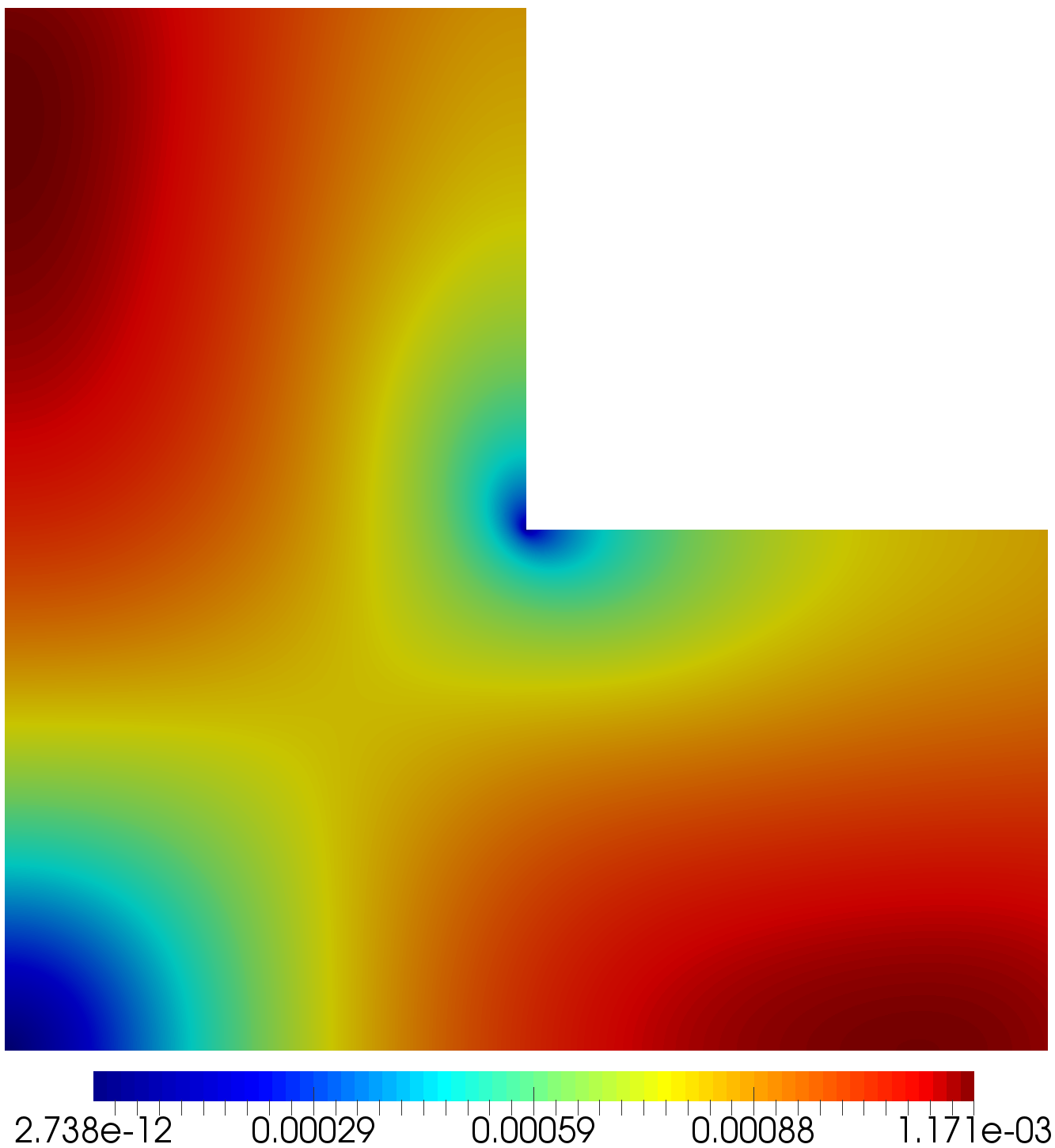}
	}
    \qquad 
	\subcaptionbox{Pressure $p$ at time t = 3.0. \label{fig:L_solution}}
	[0.19\columnwidth]
	{\includegraphics[width=0.17\textwidth,keepaspectratio]{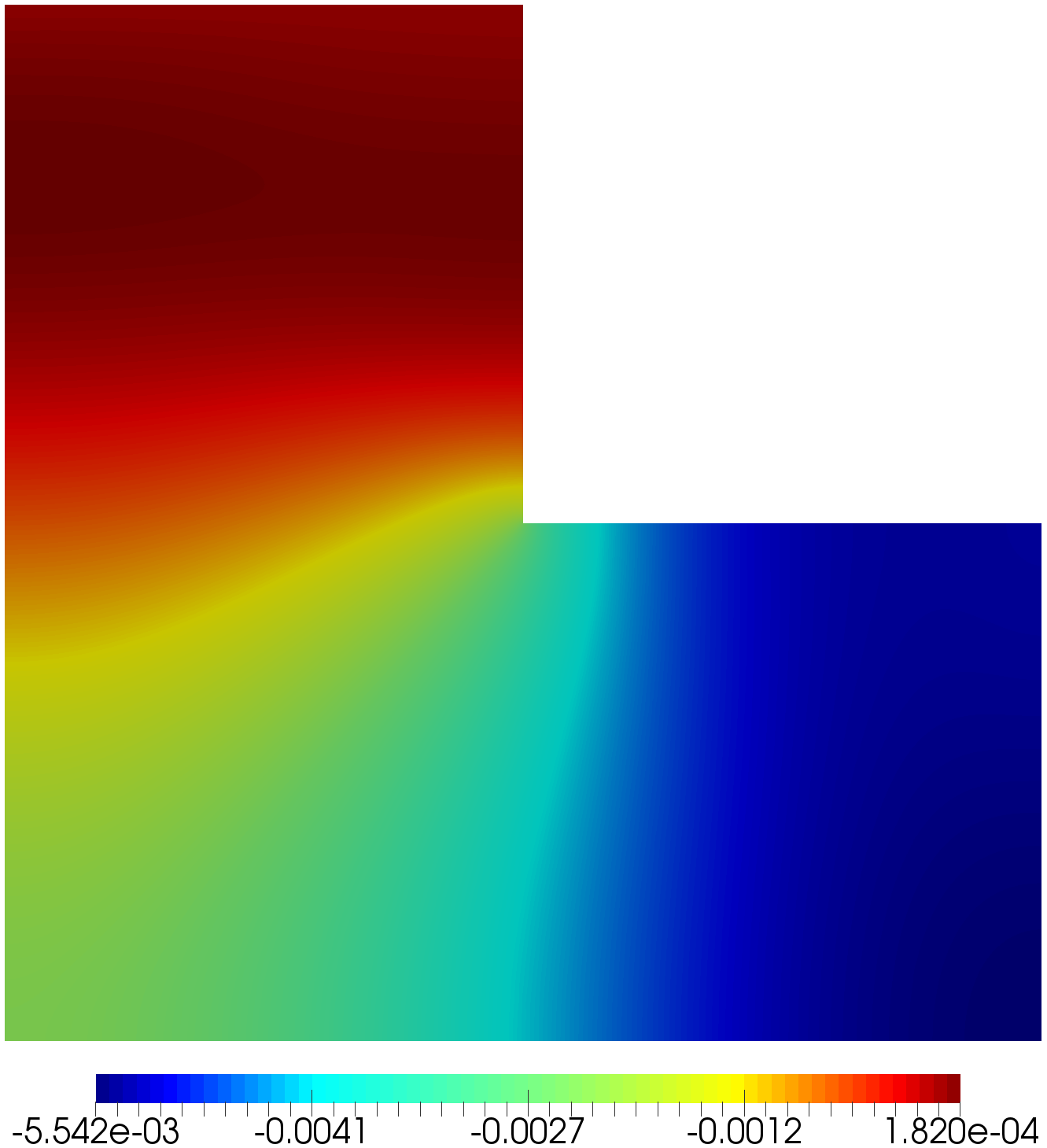}
	}
	\caption{Problem setting with boundary conditions for $\vec u$ and profile of the solution at time t = 3.8.}
	\label{fig:L_shaped_domain}
\end{figure}

Fig.~\ref{fig:ST-Conv_DG-CG} and \ref{fig:CompDGCG} illustrate and compare the results of our computations for the schemes presented in Problem \ref{Prob:DG} and \ref{Prob:CG}, respectively.  In Fig.~\ref{fig:ST-Conv_DG-CG}, the convergence in space and time of the goal quantity $G_{\vec u}$ of  \eqref{Eq:DefGQ} is clearly observed for either families of discretizations. In Fig.~\ref{fig:CompDGCG}, the superiority of the higher order members of the dG($k$) and cG($k$) schemes is illustrated. The scheme dG(1)  and, on the coarser time mesh, the scheme cG(2) are strongly erroneous. For the poroelasticity system \eqref{Eq:HPS}, this illustrates the sensitivity of numerical predictions with respect to the time discretization and argues for higher order approaches that, however, put an additional facet of  complexity on the efficient (iterative) solution of the algebraic equations; cf.~\cite{ABMS23}. Finally, in Tab.~\ref{Tab:CharNums} characteristics of the computed goal quantities are summarized.

\begin{figure}[h!]
	\centering
	\includegraphics[width=0.9\textwidth,keepaspectratio]{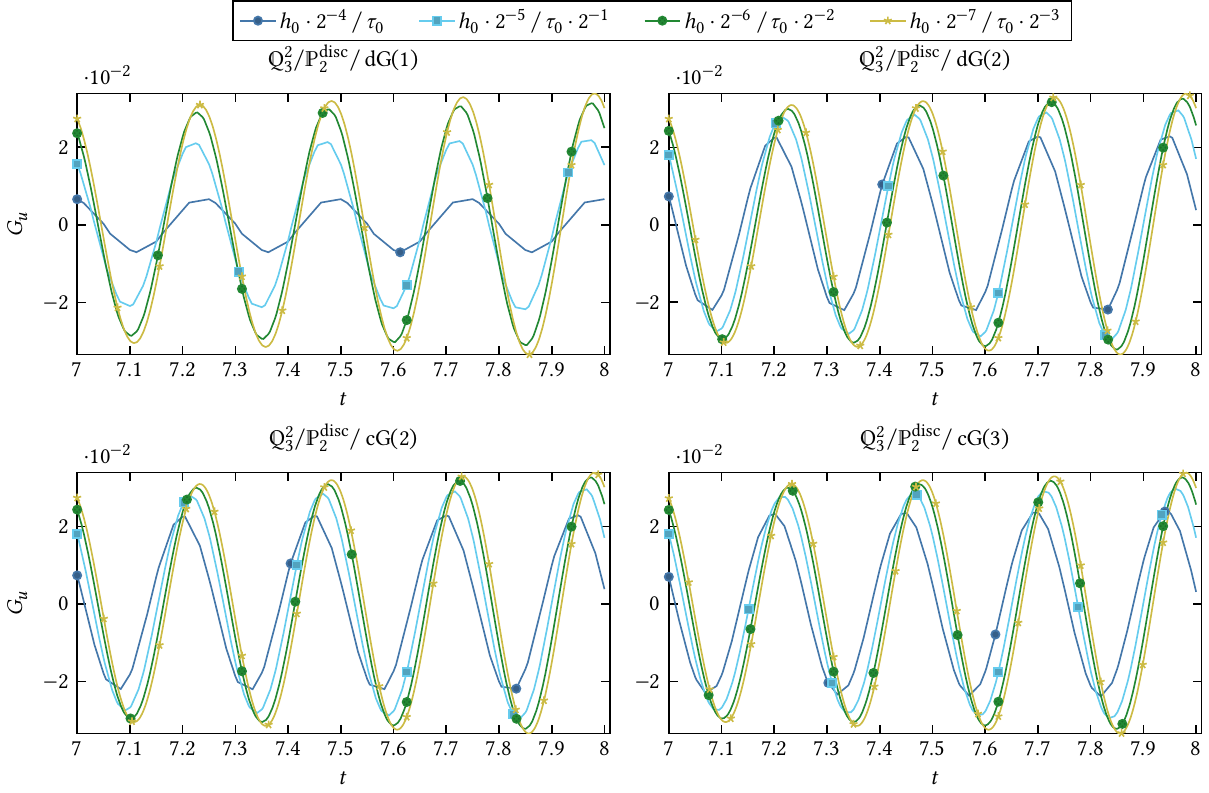}
	\caption{Computed goal quantity $G_{\vec u}$ of \eqref{Eq:DefGQ}: Space-time convergence and accuracy of the schemes defined in Problem \ref{Prob:DG} and \ref{Prob:CG}.}
	\label{fig:ST-Conv_DG-CG}
\end{figure}

\begin{figure}[h!]
	\centering
	\includegraphics[width=0.9\textwidth,keepaspectratio]{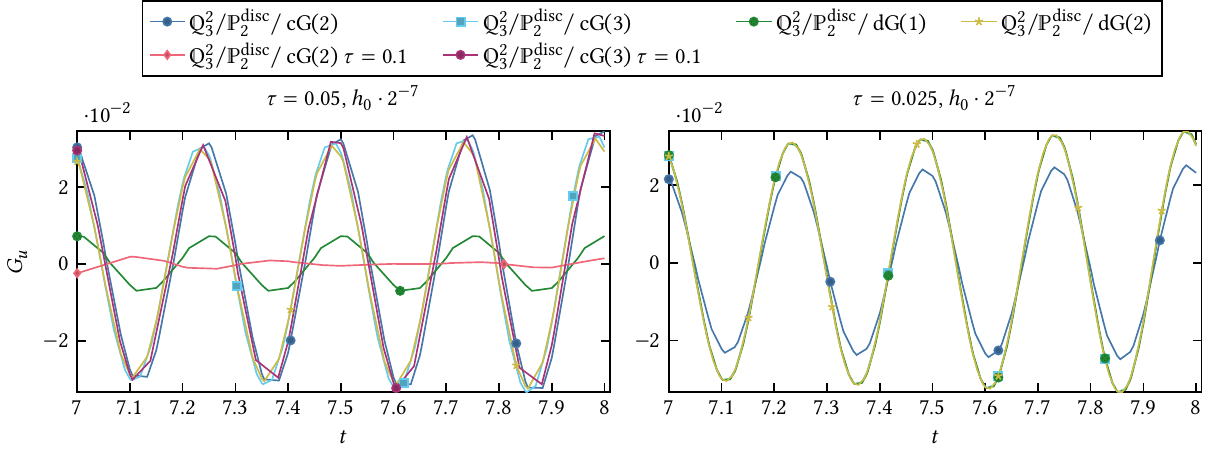}
	\caption{Computed goal quantity $G_{\vec u}$ of  \eqref{Eq:DefGQ}: Comparison of dG(1) and dG(2) of Problem \ref{Prob:DG} with cG(2) and cG(3) of Problem \ref{Prob:CG}.}
	\label{fig:CompDGCG}
\end{figure}

\begin{table}[h!]
	\centering%
	\let\mc\multicolumn%
	\setlength{\tabcolsep}{6pt}%
	\begin{tabular}{@{}cc@{\hskip 3pt}cc@{\hskip 3pt}cc@{\hskip 3pt}cc@{\hskip
				3pt}cc@{\hskip 3pt}cc@{\hskip 3pt}cc@{\hskip 3pt}cc@{\hskip 3pt}c@{}}
		\toprule
		&  \mc{2}{c}{\(\min\{G_p(t) \mid t\in [7,8]\} \)} & \mc{2}{c}{\(\max\{G_p(t) \mid t\in [7,8] \}\)} & \mc{2}{c}{\(\min\{G_{\vec u}(t)\mid t\in [7,8]\} \)} & \mc{2}{c}{\(\max\{G_{\vec u}(t)\mid t\in [7,8] \} \)}\\
		\cmidrule(lr{.66em}){2-3}\cmidrule(lr{.66em}){4-5}\cmidrule(lr{.66em}){6-7}\cmidrule(lr{.66em}){8-9}
		\(h\hphantom{_{7}}\;\;\tau\hphantom{_{3}}\) & {$cG(3)$} & {$dG(2)$} & {$cG(3)$} & {$dG(2)$} & {$cG(3)$} & {$dG(2)$} & {$cG(3)$} & {$dG(2)$}\\
		\cmidrule(r{.66em}){1-1}\cmidrule(lr{.66em}){2-2}\cmidrule(lr{.66em}){3-3}\cmidrule(lr{.66em}){4-4}\cmidrule(lr{.66em}){5-5}\cmidrule(lr{.66em}){6-6}\cmidrule(lr{.66em}){7-7}\cmidrule(lr{.66em}){8-8}\cmidrule(lr{.66em}){9-9}
		\(h_4\;\;\tau_0\) & {-2.030e-02} & {-1.903e-02} & {2.046e-02} & {1.938e-02} & {-2.387e-02} & {-2.204e-02} & {2.386e-02} & {2.290e-02}\\[0pt]
		\(h_5\;\;\tau_1\) & {-2.581e-02} & {-2.534e-02} & {2.606e-02} & {2.570e-02} & {-2.947e-02} & {-2.937e-02} & {2.981e-02} & {2.964e-02}\\[0pt]
		\(h_6\;\;\tau_2\) & {-2.850e-02} & {-2.848e-02} & {2.876e-02} & {2.882e-02} & {-3.222e-02} & {-3.222e-02} & {3.266e-02} & {3.264e-02}\\[0pt]
		\(h_7\;\;\tau_3\) & {-2.984e-02} & {-2.985e-02} & {3.024e-02} & {3.025e-02} & {-3.345e-02} & {-3.350e-02} & {3.393e-02} & {3.399e-02}\\[0pt]
		\bottomrule
	\end{tabular}
	\caption{Convergence and comparison of characteristics of the goal quantities $G_p$ and $G_{\vec u}$ in \eqref{Eq:DefGQ} for dG(2)  and cG(3) with  $\mathbb Q_3^2/\mathbb P_{2}^{\text{disc}}$ space discretization and 
		\(h_j=h_0\cdot 2^{-j}\) and \(\tau_j=\tau_0\cdot 2^{-j}\).}  
	\label{Tab:CharNums}
\end{table}

\clearpage 

\section*{Acknowledgements}
Computational resources (HPC-cluster HSUper) have been provided by the project hpc.bw, funded by dtec.bw --- Digitalization and Technology Research Center of the Bundeswehr. dtec.bw is funded by the European Union --- NextGenerationEU.

\vspace{\baselineskip}

\end{document}